\documentclass[conference]{IEEEtran}
\IEEEoverridecommandlockouts
\usepackage{cite}
\usepackage{amsmath,amssymb,amsfonts, mathtools}
\usepackage{algorithmic}
\usepackage{graphicx}
\usepackage{textcomp}
\usepackage{xcolor}
\usepackage{enumitem}
\def\BibTeX{{\rm B\kern-.05em{\sc i\kern-.025em b}\kern-.08em
    T\kern-.1667em\lower.7ex\hbox{E}\kern-.125emX}}

\newtheorem{proposition}{Proposition}[section]
\newtheorem{theorem}{Theorem}[section]

\def\weak{\rightharpoonup}

\begin{document}

\title{A Dirichlet inclusion problem on Finsler manifolds}

\author{\IEEEauthorblockN{\'Agnes Mester}
\IEEEauthorblockA{\textit{Faculty of Mathematics and Computer Science} \\
\textit{Babe\c{s}-Bolyai University}\\
Cluj-Napoca, Romania \\
\textit{\& John von Neumann Faculty of Informatics} \\
\textit{\'Obuda University}\\
Budapest, Hungary \\
agnes.mester@ubbcluj.ro}
\and
\IEEEauthorblockN{K\'aroly Szil\'ak}
\IEEEauthorblockA{\textit{Doctoral School of Applied Informatics} \\ \textit{and Applied Mathematics} \\
\textit{\'Obuda University}\\
Budapest, Hungary \\
karoly.szilak@stud.uni-obuda.hu}
}

\maketitle

\begin{abstract}
 \noindent In this paper we study a Dirichlet-type differential inclusion involving the Finsler-Laplace operator on a complete Finsler manifold. Depending on the positive $\lambda$ parameter of the inclusion, we establish non-existence, as well as existence and multiplicity results by applying non-smooth variational methods. The main difficulties are given by the problem's highly nonlinear nature due to the general Finslerian setting, as well as the non-smooth context.  
\end{abstract}


\section{Introduction and main results}

Numerous geometric and physical phenomena can be addressed by examining certain elliptic PDEs. 
An important class of such problems is represented by the elliptic equation
\begin{equation*}
    Lu(x) =  g(u(x)), \ x \in \Omega, 
    \eqno{({\mathcal P_1})} 
\end{equation*}
where $\Omega \subset \mathbb R^n$ denotes an open set $(n\geq 2)$, $L$ represents an elliptic operator and $g: \mathbb R \to \mathbb R$ is a nonlinear function verifying certain regularity and growth conditions. 
% For example, Ambrosetti and Rabinowitz \cite{AR} consider the model Dirichlet problem
% \begin{equation*}   
%     \left\{ \begin{array}{lll}
% 	-\Delta u(x)= g(u(x)),& x \in \Omega, \\
% 	u\geq 0,& x \in \Omega, \\
% 	u= 0,& x \in \partial \Omega, 
%     \end{array}\right. \eqno{({\mathcal P_2})}
% \end{equation*}   
% and guarantee the existence of one nontrivial solution by 
% applying variational arguments, where $\Delta$ is the Laplace operator, $\Omega \subset {\mathbb R}^n$ denotes a bounded open domain $(n\geq 2)$, and $g:\mathbb R\to \mathbb R$ is a continuous function with certain growth conditions at the origin and infinity. 

Due to the recent advances in geometric analysis, several elliptic problems have been studied on non-Euclidean spaces such as Riemannian and Finsler manifolds, see e.g., Bonanno, Bisci and Rădulescu \cite{BonannoBisciRadulescu}, Cheng and Yau \cite{Cheng}, Farina, Sire and Valdinoci \cite{Farina},
Farkas, Krist\'aly and Varga \cite{FarkasKristalyVarga}, Kristály, Mester and Mezei \cite{KristalyMester}, Kristály and Rudas\cite{KristalyRudas}, and references therein. 

Concerning physical applications, a significant extension of problem  $({\mathcal P_1})$ addresses the case when the nonlinear function $g$ is not continuous. This setting is handled by the replacement of $({\mathcal P_1})$ with the differential inclusion problem  
\begin{equation*}
    Lu(x) \in  \partial G(u(x)), \ x \in \Omega,
    \eqno{({\mathcal P_2})} 
\end{equation*} 
where $G:\mathbb{R} \to \mathbb{R}$ is a locally Lipschitz function and $\partial G$ stands for the generalized gradient of $G$, see e.g., Carl and Le \cite{CarlLe}, Kristály, Mezei and Szilák \cite{oscillatory, KMSZ}, and references therein.  

Motivated by these facts, 
%and inspired by the model problem $({\mathcal P_2})$,
we consider the following Dirichlet-type differential inclusion problem on a complete Finsler manifold $(M,F)$, namely
\begin{equation*} 
    \left\{ \begin{array}{lll}
        -\Delta_F u(x) \in \lambda \partial G(u(x)),& x \in \Omega, \\
        u(x) = 0, & x \in \partial\Omega,
    \end{array}\right. \eqno{({\mathcal DI})_\lambda}    
\end{equation*}
where $\Omega \subset M$ is an open, bounded set with sufficiently smooth boundary, 
$\Delta_F$ denotes the Finsler-Laplace operator defined on $(M,F)$ and $\lambda > 0$ is a parameter. Furthermore, $G:\mathbb{R} \to \mathbb{R}$ is a locally Lipschitz function and $\partial G$ represents the generalized gradient of $G$ in the sense of Clarke \cite{Clarke}.

In this case, an element $u \in W^{1,2}_{0,F}(\Omega)$ is said to be a weak solution of problem $({\mathcal DI})_\lambda$ if there exists a measurable mapping $x \mapsto \xi_x \in \partial G(u(x))$ on $\Omega$ such that for every test function  $\varphi \in C_0^\infty(\Omega)$, the correspondence  $x \mapsto \xi_x \varphi(x)$ belongs to $L^1(\Omega)$ and  
\begin{equation}\label{weak_solution}
    \int_\Omega (D\varphi(x))\big(\nabla_F u(x)\big)\,\mathrm{d}v_F(x) = \lambda \int_\Omega \xi_x \varphi(x)\,\mathrm{d}v_F(x).
\end{equation}

Considering the nonlinear term $G$, we further require the following conditions:
\begin{enumerate}[leftmargin=*,labelindent=1.7em]
    \item[$(\textbf{H}_1):$]\label{A1} $\displaystyle\lim_{t\to 0}\frac{\max\{|\xi|: \xi \in \partial G(t)\}}{t}=0$;
    \item[$(\textbf{H}_2):$]\label{A2} $\displaystyle\lim_{t\to \infty}\frac{\max\{|\xi|: \xi \in \partial G(t)\}}{t}=0$;
    \item[$(\textbf{H}_3):$]\label{A3} $\displaystyle G(0) = 0$.
\end{enumerate}

The purpose of this paper is to prove non-existence and existence/multiplicity results for the inclusion problem $({\mathcal DI})_\lambda$ depending on the parameter $\lambda$. The ambient space is given by a complete Finsler manifold, on which we consider an open, bounded subset $\Omega$. Note that the boundedness of $\Omega$ guarantees the validity of certain continuous and compact Sobolev embeddings, as well as the fact that the uniformity constant on $\Omega$, $l_{F,\Omega}$ is nonzero. All of these are crucial in our arguments. Further ingredients of the proof involve variational methods, in particular the non-smooth mountain pass theorem.

Our result can be stated as follows: 

\begin{theorem}\label{theorem1}
    Let $(M,F)$ be an $n$-dimensional complete Finsler manifold with $n \geq 3$, and let $\Omega \subset M$ be an open, bounded set with sufficiently smooth boundary. Let us consider the parameter-dependent Dirichlet inclusion problem $({\mathcal DI})_\lambda$, where $\lambda > 0$  and $G$ is a locally Lipschitz function verifying assumptions $(\textbf{H}_1) - (\textbf{H}_3)$. 
    
    Then there exist the numbers $0 < \lambda_1 < \lambda_2$ such that    
    \begin{enumerate}[leftmargin=*,labelindent=1.2em]
        \item[$(\textbf{R}_1)$] (Non-existence): $({\mathcal DI})_\lambda$ has only the trivial solution if $0 < \lambda < \lambda_1$;
        \item[$(\textbf{R}_2)$] (Existence/multiplicity): $({\mathcal DI})_\lambda$ has two different nontrivial nonnegative weak solutions when $\lambda > \lambda_2$.
    \end{enumerate}
\end{theorem}

The organization of the paper is the following. Section \ref{Sec2} provides a summary of preliminary notions and results concerning Finsler manifolds. 
In Section \ref{Sec3}, we recall some fundamental definitions and theorems of non-smooth analysis that find application in our proof. Finally, section \ref{Sec4} contains the proof of Theorem \ref{theorem1}.

\section{Preliminaries on Finsler manifolds} \label{Sec2}

In this section we briefly present the fundamental notions of Finsler manifolds, which are necessary for our further developments. For more details we refer to Bao, Chern and Shen \cite{BCS},  Ohta and Sturm \cite{OhtaSturm}, and Shen \cite{Shen01}.

Let $M$ be a connected $n$-dimensional differentiable manifold. The collection of vectors tangent to $M$ is denoted as the tangent bundle of $M$, defined as
$$TM=\cup_{x \in M}\{(x,v): v \in T_{x} M\},$$ 
where $T_{x} M$ represents the tangent space of $M$ at the point $x$.

The pair $(M,F)$ is called a Finsler manifold if $F: TM \to [0,\infty)$ is a continuous function verifying the following conditions:
\begin{enumerate}[label=(\roman*)]
	\item $F \in C^{\infty}(TM \setminus \{ 0 \})$; 
	\item $F(x,\lambda v) = \lambda F(x,v)$, for every $\lambda \geq 0$ and $(x,v) \in TM$;
	\item the Hessian matrix $\left[\left(  \frac{1}{2}F^{2}(x,v)\right) _{v^{i}v^{j}} \right]_{i,j=\overline{1,n}}$
	is positive definite for every $(x,v)\in TM \setminus \{0\}.$
\end{enumerate}
In this case we say that $F$ is a Finsler metric on $M$.

If, in addition, $F(x,\lambda v) = |\lambda| F(x,v)$ holds for all $\lambda \in \mathbb{R}$ and $(x,v) \in TM$, then  $F$ is said to be symmetric and the Finsler manifold is called reversible. Otherwise, $(M,F)$ is classified as nonreversible.

The reversibility constant of $(M,F)$ is defined by the number 
\begin{equation*} 
r_{F} = \sup_{x\in M} ~ \sup_{\substack{ v\in T_{x} M \setminus \{0\}}} \frac{F(x,v)}{F(x,-v)} ~ \in [1, \infty],
\end{equation*}
and it measures how much the Finsler manifold deviates from being reversible, see Rademacher \cite{Rademacher}.
In particular, $r_F = 1$ if and only if $(M,F)$ is a reversible Finsler manifold. 

The uniformity constant of $(M,F)$ is defined by the number
\begin{equation*}
 l_{F} = \inf_{x\in M} ~ \inf_{v,w,z\in  T_x M\setminus \{0\}}\frac{g_{(x,w)}(v,v)}{g_{(x,z)}(v,v)} ~ \in ~ [0, 1],
\end{equation*}
which measures how much $F$ deviates from being a Riemannian structure, see Egloff \cite{Egloff}. Here $g$ denotes the fundamental tensor of $(M,F)$, see Bao, Chern and Shen \cite{BCS}.  We have $l_F = 1$ if and only if $(M,F)$ is a Riemannian manifold, see Ohta \cite{Ohta}.

A $C^\infty$-differentiable curve $\gamma: [a,b] \to M$ is called a geodesic if its velocity field $\dot \gamma$ is parallel along the curve.
The Finsler manifold is said to be complete if every geodesic segment $\gamma: [a,b] \to M$ can be extended to a geodesic defined on $\mathbb{R}$. 

The dual metric $F^*:T^*M \to [0,\infty)$ is called the polar transform of $F$, defined as
\begin{equation*}  
F^*(x,\alpha) = \sup_{v \in T_xM \setminus \{0\}} ~ \frac{\alpha(v)}{F(x,v)},
\end{equation*}
where $T^*M = \bigcup_{x \in M}T^*_{x} M $ is the cotangent bundle of $M$ and $T^*_{x} M$ represents the dual space of $T_{x} M$.

In local coordinates, the Legendre transform $J^*:T^*M \to TM$ is defined as
\begin{equation*}  
J^*(x,\alpha) = \sum_{i=1}^n \frac{\partial}{\partial \alpha_i}\left(\frac{1}{2} F^{*2}(x,\alpha)\right)\frac{\partial}{\partial x^i}.
\end{equation*}
Note that if $ J^*(x, \alpha) = (x,v) $, then one has 
\begin{equation*} 
F(x,v) = F^*(x,\alpha) \quad \text{and} \quad \alpha(v) = F^*(x,\alpha) F(x,v).
\end{equation*}
%One has $F(J^*(x,\alpha)) = F^*(x,\alpha)$.

Furthermore, if $u: M \to \mathbb R$ is a function of class $C^1$, then the gradient of $u$ is defined as 
\begin{equation*}  
\nabla_F u(x) = J^*(x, Du(x)), \ \forall x \in M,
\end{equation*}
where $Du(x) \in T_x^*M$ denotes the differential of $u$ at the point $x$. In particular, one has that
$$F(x, \nabla_F u(x)) = F^*(x, Du(x)), \ \forall x \in M.$$

It is important to note that, in general, $\nabla_F$ is a nonlinear operator. 
Also, the mean value theorem implies that
\begin{align}\label{mean_theorem}
    &(Du(x) - Df(x))(\nabla_F u(x)-\nabla_F f(x)) \geq \nonumber \\
    & \geq l_F F^{*2}(Du(x) - Df(x)),
\end{align}
for every $u,f \in C^1(M)$.

For a function $u \in C^2(M)$, the Finsler-Laplace operator  
$\Delta_F$ is defined as 
$$\Delta_F u(x) = \mathrm{div}(\nabla_F u(x)),$$
where 
\begin{equation*} 
\mathrm{div} V(x) = 
\frac{1}{\sigma_F(x)} \sum_{i=1}^n \frac{\partial}{\partial x^i} \Big(\sigma_F(x) V^i(x) \Big)
\end{equation*}
for some vector field $V$ on $M$. Here, $\sigma_F$ represents the density function defined by
$\sigma_F(x) = \frac{ \omega_n}{\mathrm{Vol}(B_x(1))},$
where $\omega_n$ and $\mathrm{Vol}(B_x(1))$ denote the Euclidean volume of the $n$-dimensional Euclidean unit ball and the set 
$$B_x(1) = \Big\{ (v^i) \in \mathbb{R}^n :~ F \Big(x, \sum_{i=1}^n v^i \frac{\partial}{\partial x^i} \Big) < 1 \Big\} \subset \mathbb{R}^n,$$
respectively.
Again, the Finsler-Laplace operator $\Delta_F$ is generally nonlinear.

The canonical Hausdorff volume form $\mathrm{d}v_F$ on $(M, F)$ is defined as
\begin{equation*} \label{Hausdorff_measure}
\mathrm{d}v_F(x) = \sigma_F(x) \mathrm{d} x^1 \land \dots \land \mathrm{d} x^n.
\end{equation*}
In the following, we may omit the parameter $x$ for the sake of brevity.
The Finslerian volume of an open set $\Omega\subset M$ is given by $\mathrm{Vol}_F(\Omega) = \int_\Omega \mathrm{d}v_F(x)$.

Let $\Omega \subset M$ be an open set.
The Sobolev space $W^{1,2}_F(\Omega)$ associated with the Finsler structure $F$ and the canonical Hausdorff measure $\mathrm{d}v_F$ is defined as
$$W^{1,2}_F(\Omega) = \left\{u\in W^{1,2}_\mathrm{loc}(\Omega): \int_\Omega {F^*}^2(x,Du(x))\mathrm{d}v_F < \infty \right\}.$$

It can be proved that $W^{1,2}_F(\Omega)$ is the closure of $C^\infty(\Omega)$ with respect to the (generally asymmetric) norm 
\begin{equation*}\label{Sobolev_norm1}
\|u\|_{W^{1,2}_F(\Omega)}=\left(\int_\Omega {F^*}^2(x,Du(x))\,\mathrm{d}v_F + \int_\Omega |u(x)|^2\,\mathrm{d}v_F \right)^{\frac{1}{2}}.
\end{equation*}
% which is also equivalent with the norm 
% \begin{equation*}\label{Sobolev_norm2}
% \left(\int_\Omega {F^*}^2(x,Du(x))\,\mathrm{d}v_F(x)\right)^{\frac{1}{2}} + \left(\int_\Omega |u(x)|^2\,\mathrm{d}v_F(x) \right)^{\frac{1}{2}}.
% \end{equation*}

The space $W_{0,F}^{1,2}(\Omega)$ is defined as the closure of $C_{0}^{\infty}(\Omega)$ with respect to the norm $\|\cdot\|_{W^{1,2}_F(\Omega)}$.

Since $F$ is not necessarily symmetric, the Sobolev spaces $W^{1,2}_F(\Omega)$ and $W^{1,2}_{0,F}(\Omega)$ are generally asymmetric normed spaces. 
However, if $(M, F)$ is a complete Finsler manifold and $\Omega$ is a bounded subset of $M$, then due to Farkas, Krist\'aly and Varga \cite{FarkasKristalyVarga}, we have that $W^{1,2}_F(\Omega)$ and $W^{1,2}_{0,F}(\Omega)$ are reflexive biBanach spaces (i.e., complete asymmetric normed spaces, see Cobza\c s \cite{Cobzas}).

Furthermore, when $\Omega \subset M$ is an open, bounded set with sufficiently smooth boundary and $n = \text{dim} M \geq 3$, the Sobolev space $W^{1,2}_{F}(\Omega)$  can be continuously embedded into the Lebesgue space $L^q(\Omega)$ (see Hebey \cite[Theorem 2.6]{Hebey}) for every $q \in [1,2^*]$, where $2^* = \frac{2n}{n-2}$. In particular, for every $q \in [1,2^*]$, there exists a constant $C_q > 0$ such that
\begin{equation} \label{cont_embedding} 
\left(\int_\Omega |u(x)|^q\,\mathrm{d}v_F\right)^\frac{1}{q} \leq C_q \left(\int_\Omega F^{*2}(x,Du(x))\,\mathrm{d}v_F\right)^\frac{1}{2},
\end{equation}
for every $u \in W^{1,2}_{F}(\Omega)$. Moreover, %under similar assumptions on $\Omega$ and $n$, 
the compact embeddings $W^{1,2}_{F}(\Omega) \hookrightarrow L^q(\Omega)$ are also valid for every $q \in [1,2^*)$, see Hebey \cite[Theorem 2.9]{Hebey}. 

The operators $\mathrm{div}$ and $\Delta_F$ can be defined in a distributional sense as well, see Ohta and Sturm \cite{OhtaSturm}. For instance, for a function $u \in W^{1,2}_F(\Omega)$, $\Delta_F u$ is defined as
\begin{equation*}  \label{divergence_theorem}
\int_{\Omega} \varphi(x) \Delta_F u(x) ~ \mathrm{d}v_F = -\int_{\Omega} D\varphi(x)\big( \nabla_F u(x)\big) \mathrm{d}v_F, 
\end{equation*}
for all $\varphi \in C^{\infty}_0(\Omega)$.

\section{Elements of non-smooth analysis} \label{Sec3}

In this segment, we review the fundamental characteristics of locally Lipschitz functions that find application in our proofs. For more details, see Clarke \cite{Clarke}.

Let $X$ be a real Banach space equipped with the norm $\|\cdot\|$.  A function $h:X \rightarrow\mathbb{R}$ is said to be locally Lipschitz if every point $u\in X$ has a neighborhood $N_u \subset X$ such that
\begin{equation*}%\label{loc-lip}
    \vert h(u_{1})-h(u_{2})\vert\leq K\| u_{1}-u_{2}\|, \quad \forall u_{1}, u_{2} \in N_u,
\end{equation*}
where $K > 0$ is a constant depending on $N_u$. 

Now suppose that $h:X\rightarrow\mathbb{R}$ is a locally Lipschitz function. 
The generalized directional derivative of $h$ at  $u\in X$ in the direction $v\in X$ is defined as
$$ h^{0}(u;v):=\limsup_{\scriptstyle {\it w\rightarrow u}\atop 	\scriptstyle \it t\searrow 0}\frac{h(w+tv)-h(w)}{t}.$$
Note that if $h$ is of class $C^1$ on $X$, then $ h^{0}(u;v) = \langle h'(u), v\rangle$ for all $u, v \in X$.
Hereafter, $\langle \cdot, \cdot \rangle$ and $\|\cdot\|_*$ denote the duality mapping on $(X^*, X)$ and the norm of the dual space $X^*$. 

The Clarke subdifferential of the locally
Lipschitz function $h:X\rightarrow\mathbb{R}$ at a
point $u\in X$  is defined by the set
$$\partial h(u):=\left\{\zeta\in X^{\ast}: \langle \zeta, v\rangle\leq h^{0}(u; v),\ \forall v\in X \right\}.$$
	
An element	$u\in X$ is called a critical point of $h$ if $0\in \partial h(u)$, see Chang \cite[Definition 2.1]{Chang}. 

For a locally Lipschitz function, the following assertions are available.

\begin{proposition}\label{prop-lok-lip-0} (Clarke \cite{Clarke}) \it 
    Let $h: X\rightarrow\mathbb{R}$ be a locally Lipschitz function. Then the following properties hold:
    \begin{enumerate} 
        \item $($Lebourg's mean value theorem$)$ Let  $U$  be an open subset of a Banach space $X$ and $u, v$ be
    two points of $U$ such that the line          segment $[u,v] = \{(1-t)u+tv: t \in     [0,1]\}$ is in $U$. If  $h:U\rightarrow\mathbb{R}$ is a Lipschitz function, then there exist a point $w\in (u,v)$ and $\zeta\in \partial h(w)$ such that 
        $h(v)-h(u)\in \langle \zeta,v-u\rangle.$ 
    
        \item If $j:X\to \mathbb R$, $j \in C^1(X)$, then $\partial (j+h)(u)=j'(u)+\partial h(u)$ and 
        $(j+h)^0(u;v)=\langle j'(u),v\rangle +h^0(u;v)$ for every $u,v\in X.$ 
        \item $(-h)^0(u;v)=h^0(u;-v)$ for every $u,v\in X.$ 
        \item $\partial (\alpha h)(u)=\alpha \partial h(u)$ for every $\alpha\in \mathbb R$ and $u\in X.$
        
        %			\item The set-valued map $\partial h:X\rightsquigarrow {X^{\ast}}$
        %			is u.s.c. from $s-X$ into $w^\ast-X^\ast$.
    \end{enumerate}
\end{proposition}

\section{Proof of Theorem \ref{theorem1}} \label{Sec4}

    \subsection{Proof of the non-existence result $(\textbf{R}_1)$}

    Let $u \in W^{1,2}_{0,F}(\Omega)$ be a weak solution of $(\mathcal DI)_\lambda$. By density arguments, we can choose $\varphi \coloneqq u$ in \eqref{weak_solution}, thus we obtain
    \begin{equation}\label{egyik}
    \int_\Omega {F^*}^2(x,Du(x))\,\mathrm{d}v_F = \lambda \int_\Omega \xi_x u(x)\,\mathrm{d}v_F,
    \end{equation}
    where $\xi_x \in \partial G(u(x))$ for every $x \in \Omega$ such that  $x \mapsto \xi_x u(x)$ belongs to $L^1(\Omega)$. 
    
    Based on the assumptions $(\textbf{H}_1)$ and $(\textbf{H}_2)$, for every $\varepsilon>0$, one can find the numbers $\delta_1, \delta_2>0$ such that
    \begin{equation} \label{A1_and_A2_concl} 
        |\xi|\leq \varepsilon t, ~ \forall \xi \in \partial G(t) \text{ and } \forall  0 < t < \delta_1 \textnormal{ or } t  > \delta_2.  
    \end{equation} 

    Then for every $\varepsilon>0$ we have that 
    \begin{equation}\label{masik}
    \int_\Omega \xi_x u(x)\,\mathrm{d}v_F \leq  \varepsilon \int_\Omega |u(x)|^2\,\mathrm{d}v_F.
    \end{equation}
    
    Consequently, from \eqref{egyik}, \eqref{masik}  and the continuous embedding $W^{1,2}_{F}(\Omega) \subset L^2(\Omega)$ (see  \eqref{cont_embedding}), it follows that  
    \begin{equation*}
    (\lambda \varepsilon C_2^2 -1 )\int_\Omega {F^*}^2(x,Du(x))\,\mathrm{d}v_F \geq  0,
    \end{equation*}
    where $C_2$ denotes the embedding constant from \eqref{cont_embedding} in the case $q=2$.

    Accordingly, if $\lambda \varepsilon C_2^2 < 1$, i.e., when
    $$\lambda < \frac{1}{\varepsilon C_2^2} \eqqcolon \lambda_1,$$
    then we necessarily have $u = 0$ a.e. on $\Omega$, which concludes the proof of $(\textbf{R}_1)$.

    \subsection{Proof of the existence/multiplicity result $(\textbf{R}_2)$} \label{B}

    Having in mind the inclusion $(\mathcal DI)_\lambda $,   let us construct the modified problem 
    \begin{equation*} 
    \left\{ \begin{array}{lll}
        -\Delta_F u(x) \in \lambda \partial \widetilde{G}(u(x)),& x \in \Omega, \\
        u(x) = 0, & x \in \partial\Omega,
    \end{array}\right. \eqno{(\widetilde{\mathcal DI})_\lambda}    
    \end{equation*}
    where we define $\widetilde{G}: \mathbb R \to \mathbb R$ as 
    $$\widetilde{G}(t) =  \begin{cases}
    0, & \text{ if } t < 0 \\
    G(t), & \text{ if } t \geq 0. 
    \end{cases}$$
    Clearly, $\widetilde{G}$ is a locally Lipschitz function which also satisfies hypotheses $(\textbf{H}_1)-(\textbf{H}_3)$.

    We consider the energy functional associated with problem $(\widetilde{\mathcal DI})_\lambda$ for every $\lambda >0$, namely
    $$\mathcal{E}_\lambda: W^{1,2}_{0,F}(\Omega)\to \mathbb{R}, \quad \mathcal{E}_\lambda(u) = \frac{1}{2} \mathcal{N}(u) - \lambda \mathcal{G}(u),$$ 
    where $\mathcal{N}, \mathcal{G}: W^{1,2}_{0,F}(\Omega)\to \mathbb{R}$,
    $$\mathcal{N}(u) = \int_\Omega{F^*}^2(x,Du(x))\,\mathrm{d}v_F $$ and
    $$\mathcal{G}(u)=\int_\Omega \widetilde{G}(u(x))\,\mathrm{d}v_F.$$

    Let us observe that if $u \in W^{1,2}_{0,F}(\Omega)$ is a weak solution of $(\widetilde{\mathcal DI})_\lambda$, then $u$ is a nonnegative weak solution of the original problem $(\mathcal DI)_\lambda$. Indeed, suppose that there exists $u \in W^{1,2}_{0,F}(\Omega)$ a nontrivial solution of $(\widetilde{\mathcal DI})_\lambda$ for some $\lambda$. 
    Let $u_-(x) = \min(0,u(x))$ and $\Omega_- = \{x \in \Omega: u(x) < 0\}$. Multiplying the first equation of $(\widetilde{\mathcal DI})_\lambda$ by $u_-$ and integrating over $\Omega$, we obtain that 
    $$\int_{\Omega_-} F^{*2}(x, Du_-(x)) \mathrm{d}v_F = 0,$$  
    which in turn yields that $u_- = 0$, since $u = 0$ on $\partial \Omega$. Hence $u \geq 0$ on $\Omega$, which means that $\widetilde{G}(u(x)) = G(u(x))$ on $\Omega$.  
    
    Therefore, it is enough to study the weak solutions of problem $(\widetilde{\mathcal DI})_\lambda$, which are given by the critical points of $\mathcal{E}_\lambda$. Accordingly, in the following subsections we shall study the properties of the energy functional $\mathcal{E}_\lambda$, namely: coercivity, boundedness from below and the validity of the non-smooth Palais-Smale condition.

    \subsection{$\mathcal{E}_{\lambda}$ is coercive and bounded below} \label{Coercivity}
    
    In this subsection we prove that the energy functional $\mathcal E_{\lambda}$ is coercive and bounded below. 

    Let $\delta >0$ be arbitrarily fixed, and let $S_1, S_2 \subset \Omega$ denote the following level sets for some function $u \in W^{1,2}_{0,F}(\Omega)$:
    \begin{align*}{}
        S_1 &= \{x \in \Omega: u(x)\leq\delta \} ~ ~ \text{ and} \\
        S_2 &= \{x \in \Omega: u(x)>\delta \}.
    \end{align*}
    
    Evidently, we have
    \begin{equation*}
        \mathcal{G}(u) = \int_{S_1} \widetilde{G}(u(x))\,\mathrm{d}v_F + \int_{S_2} \widetilde{G}(u(x))\,\mathrm{d}v_F .
    \end{equation*}
        
    For the first term, we can write
    \begin{equation*}{}
        \int_{S_1} \widetilde{G}(u(x)) {\rm d}v_F 
        \leq \max_{t \leq \delta} |\widetilde{G}(t)| \cdot \mathrm{Vol}_F(\Omega).
    \end{equation*}

    For the second term, we apply Lebourg's mean value theorem (see Proposition \ref{prop-lok-lip-0}, \textit{1)}). Then, based on relation \eqref{A1_and_A2_concl} we conclude that for every $\epsilon>0$ we can choose $\delta$ such that
    \begin{align*}
        \int_{S_2} \widetilde{G}(u(x)) {\rm d}v_F &=  \int_{S_2} \widetilde{G}(u(x))-\widetilde{G}(\delta) {\rm d}v_F  + \int_{S_2} \widetilde{G}(\delta) {\rm d}v_F \nonumber\\
        &= \int_{S_2} \langle \xi(x), u(x)-\delta\rangle {\rm d}v_F + \int_{S_2} \widetilde{G}(\delta) {\rm d}v_F \nonumber \\
        &\leq  \varepsilon \int_{\Omega} |u(x)|^2 {\rm d}v_F + \widetilde{G}(\delta)\mathrm{Vol}_F(\Omega).\nonumber
    \end{align*}
    
    Therefore, taking into account the continuous embedding $W^{1,2}_F(\Omega) \subset L^2(\Omega)$ (see \eqref{cont_embedding} in the case $q=2$), it results that  
    \begin{align*}
        \mathcal E_{\lambda}(u) &= \frac{1}{2}\mathcal{N}(u) - \lambda \mathcal{G}(u) \nonumber \\
        &\geq \left(\frac{1}{2 C_2^2} - \lambda \epsilon\right)\int_{\Omega} |u(x)|^2\,\mathrm{d}v_F \nonumber \\ 
        &- \lambda \bigg(\max_{t \leq \delta} |\widetilde{G}(t)| + \widetilde{G}(\delta)\bigg)\mathrm{Vol}_F(\Omega).
        \nonumber
    \end{align*}
    
    Since $\Omega$ is bounded, by choosing 
    $\varepsilon < \frac{1}{2C_2^2\lambda}$ sufficiently small, it follows that the energy functional $\mathcal E_{\lambda}$ is coercive and bounded from below.

    \subsection{ $\mathcal E_{\lambda}$ satisfies the non-smooth Palais-Smale condition}\label{PS}	

    Let $(u_k)_k$ be a Palais-Smale sequence for $\mathcal E_{\lambda}$ in $W^{1,2}_{0,F}(\Omega)$, i.e., suppose that $(\mathcal{E}_{\lambda}(u_k))_k$ is bounded and $m(u_k) \to 0$ as $k \to \infty$, where $$m(u) = \min \{\|\xi\|_* : \xi \in \partial \mathcal{E}_{\lambda}(u)\}.$$ 
	
    Due to the coercivity of  $\mathcal{E}_{\lambda}$, it follows that the sequence $(u_k)_k$ is bounded in $W^{1,2}_{0,F}(\Omega)$. As $W^{1,2}_{F}(\Omega)$ can be compactly embedded into $L^{q}(\Omega)$ for any $q \in [2,2^*)$ (see \eqref{cont_embedding}),
    there exists a function $u \in W^{1,2}_{0,F}(\Omega)$ such that, up to a subsequence, $u_k\weak u$ weakly in $W^{1,2}_{0,F}(\Omega)$ and $u_k \to u$ strongly in $L^q(\Omega)$ for every $q \in [2,2^*)$. 
   
    Based on Proposition \ref{prop-lok-lip-0}, we can write for the generalized directional derivative of $\mathcal E_{\lambda}$ that
    \begin{equation*}
	\mathcal E_{\lambda}^{0}(u;u_k-u) = \frac{1}{2}\langle \mathcal{N}^{'}(u), u_k-u \rangle   + \lambda (-\mathcal{G})^{0}(u;u_k-u)
    \end{equation*}
    and 
    \begin{equation*}
        \mathcal E_{\lambda}^{0}(u_k;u-u_k) = \frac{1}{2}\langle \mathcal{N}^{'}(u_k), u-u_k \rangle 
        + \lambda (-\mathcal{G})^{0}(u_k;u-u_k). 
    \end{equation*}

    Our aim is to prove that, up to a subsequence, $(u_k)_k$ strongly converges to $u$ in $W^{1,2}_{0,F}(\Omega)$. Therefore, we analyze the expression    
    \begin{align*}
	I_k &\coloneqq \int_{\Omega} (Du - Du_k)(\nabla_F u - \nabla_F u_k)\,\mathrm{d}v_F \\
%    &=  \int_{\Omega} Du(\nabla_F u - \nabla_F u_k)\,\mathrm{d}v_F \\
%    &- \int_{\Omega} Du_k(\nabla_F u - \nabla_F u_k)\,\mathrm{d}v_F \\
    &= \langle \mathcal{N}^{'}(u), u-u_k \rangle - \langle \mathcal{N}^{'}(u_k), u-u_k \rangle  \\
    &= 2 \cdot \big\{-\mathcal E_{\lambda}^{0}(u;u_k-u) + \lambda (-\mathcal{G})^{0}(u;u_k-u) \\ 
    &-\mathcal E_{\lambda}^{0}(u_k;u-u_k) + \lambda (-\mathcal{G})^{0}(u_k;u-u_k)\big\}.
    \end{align*}

    First, since $(u_k)_k$ is a Palais-Smale sequence for $\mathcal E_{\lambda}$, we have that  
    \begin{equation*}
     \mathcal E_{\lambda}^{0}(u;u_k-u) +\mathcal E_{\lambda}^{0}(u_k;u-u_k) \to 0 ~ \text{ as } ~ k\to \infty. 
    \end{equation*}

    For the remaining terms, first we recall the fact that $\partial \widetilde{G}$ is upper semicontinuous. Therefore, applying the Weierstrass theorem with conditions $(\textbf{H}_1) \& (\textbf{H}_2)$ gives us the boundedness of the function 
    \begin{equation}\label{boundedness}
    t \mapsto \frac{\max \{ |\xi|: \xi \in \partial \widetilde{G}(t) \} }{t^{p-1}}, 
    \end{equation}
    where $t >0$ and $p \in (2,2^*)$.  
	
    Combining \eqref{A1_and_A2_concl} with the boundedness of \eqref{boundedness} on $[\delta_1, \delta_2]$ yields that for every $\varepsilon>0$, 
    there exists $k_{\varepsilon} > 0$ such that
    \begin{equation} \label{estimate_1}
       |\xi| \leq \varepsilon t + k_{\varepsilon}t^{p-1}, ~ \forall t \geq 0, ~ \forall  \xi \in \partial \widetilde{G}(t).
    \end{equation}

    Accordingly, it follows that
    \begin{align*}
    S_k&\coloneqq \mathcal{G}^{0}(u;u-u_k) +  \mathcal{G}^{0}(u_k;u_k-u)  \\
    &\leq \int_\Omega \left[\widetilde{G}^{0}(u;u-u_k) + \widetilde{G}^{0}(u_k;u_k-u) \right]{\rm d}v_F \\
    &\leq  \int_\Omega \left[\max_{\xi}\{\xi(u-u_k)\} + \max_{\eta_k}\{\eta_k(u_k-u)\} \right]{\rm d}v_F \\
    &\leq  2 \varepsilon\left[\|u\|^2_{L^2(\Omega)}+\|u_k\|^2_{L^2(\Omega)}\right] \\
    &+ k_\varepsilon
    \|u-u_k\|_{L^2(\Omega)} \left[\|u\|^{p}_{L^p(\Omega)}+\|u_k\|^{p}_{L^p(\Omega)}\right]    ,
    \end{align*}
    where $\xi \in \partial \widetilde{G}(u)$ and $\eta_k \in \partial \widetilde{G}(u_k)$.
    Since $u_k \to u$ strongly in $L^q(\Omega)$ for any $q \in [2,2^*)$ it follows that $\limsup_{k \to \infty} {S_k} \leq 0$. 
    
    On the other hand, by using \eqref{mean_theorem} on $\Omega$, we obtain that
    \begin{align*}
     I_k &= \int_\Omega (Du - Du_k)(\nabla_F u-\nabla_F u_k)\,{\rm d}v_F \geq \\
    & \geq l_{F, \Omega} \int_\Omega F^{*2}(Du(x) - Du_k(x))\,{\rm d}v_F ,
    \end{align*}
    where
    \begin{equation*}
        l_{F, \Omega} = \inf_{x\in \Omega} ~ \inf_{v,w,z\in  T_x M\setminus \{0\}}\frac{g_{(x,w)}(v,v)}{g_{(x,z)}(v,v)} 
    \end{equation*}
    is the uniformity constant associated with $\Omega$. Since $\Omega \subset M$ is bounded, we can prove that $0 < l_{F, \Omega} \leq 1$, thus we obtain that
    \begin{equation*}
     \int_\Omega F^{*2}(Du(x) - Du_k(x))\,{\rm d}v_F \to 0 
    \end{equation*}
    as $k \to 0$. Therefore, $u_k \to u$ strongly in $W^{1,2}_{0,F}(\Omega)$, which proves the claim.

    \subsection{First solution}

    In this section we prove that there exists a local minimum for the energy functional $\mathcal{E}_{\lambda}$. 
	
    By applying Lebourg's mean value theorem (see Proposition \ref{prop-lok-lip-0}, \textit{1)}),  estimate \eqref{estimate_1} together with assumption $(\textbf{H}_3)$ gives the inequality
    \begin{equation*} %\label{estimate_2}
        0 \leq |\widetilde{G}(t)| \leq \varepsilon t^2 + k_{\varepsilon}|t|^p \textnormal{,   } ~ \forall t \in \mathbb{R}.
    \end{equation*} 
	
    Therefore, the continuous embeddings $W^{1,2}_{F}(\Omega) \subset L^q(\Omega)$ for every $q \in [2,2^*)$ yield that 
    \begin{align}\label{G_estimate}
        0&\leq|\mathcal{G}(u)| \leq \int_{\Omega} |\widetilde{G}(u(x))|\,{\rm d}v_{F} \nonumber\\
        &\leq \varepsilon \int_{\Omega} u(x)^2\,{\rm d}v_F  + k_{\varepsilon} \int_{\Omega} |u(x)|^p\,{\rm d}v_F  \nonumber\\
        &\leq  \varepsilon C^2_2 \int_{\Omega} F^{*2}(x,Du(x)){\rm d}v_F \\
        &+ k_{\varepsilon} C_p^p  \left(\int_{\Omega} F^{*2}(x,Du(x)){\rm d}v_F\right)^\frac{p}{2}, \nonumber
    \end{align}
    for all $u \in W^{1,2}_{0,F}(\Omega)$,
    where $C_2$ and $C_p$ denote the embedding constants from \eqref{cont_embedding} where $p \in (2, 2^*)$.
    
    Since $\varepsilon >0$ can be arbitrarily small
    and $p>2$, the previous estimate implies that for every $u \in W^{1,2}_{0,F}(\Omega) \setminus \{0\}$, we have
    \begin{equation}\label{nullaban}
        \lim_{\mathcal{N}(u)\to 0}\frac{\mathcal{G}(u)}{\mathcal{N}(u)} = 0.
    \end{equation}
    
    Similarly, it can be proven that 
    \begin{equation}\label{vegtelenben}
        \lim_{\mathcal{N}(u)\to \infty}\frac{\mathcal{G}(u)}{\mathcal{N}(u)} = 0.
    \end{equation}

    Indeed, taking into account the boundedness of the function 
    $$t \mapsto \frac{\max \{ |\xi|: \xi \in \partial \widetilde{G}(t) \} }{t^{1/2}}$$ 
    on $[\delta_1, \delta_2]$  together with relation (\ref{A1_and_A2_concl}), for every $\varepsilon >0 $ one can find an $l_{\varepsilon} > 0$ such that
    \begin{equation*} 
        |\xi| \leq \varepsilon t + l_{\varepsilon}t^{1/2}, ~ \forall t \geq 0, ~ \forall \xi \in \partial \widetilde{G}(t).
    \end{equation*}
    In a similar fashion as before, applying Lebourg's mean value theorem and the continuous Sobolev embeddings, for all $u \in W^{1,2}_{0,F}(\Omega)$ we obtain the estimate 
    \begin{align*}
        0 &\leq \mathcal{G}(u) \leq \int_{\Omega} |\widetilde{G}(u(x))| \mathrm{d}v_{F}\\
        &\leq \varepsilon \int_{\Omega} u(x)^2 \mathrm{d}v_F  + l_{\varepsilon} \int_{\Omega} |u(x)|^{\frac{3}{2}} \mathrm{d}v_F \\
        &\leq  \varepsilon C_2^2 \int_{\Omega} F^{*2}(x,Du(x)) \mathrm{d}v_F \\
        &+ l_{\varepsilon} C_{3/2}^{\frac{3}{2}} \left(\int_{\Omega} F^{*2}(x,Du(x))\mathrm{d}v_F\right)^{\frac{3}{4}},
    \end{align*}
    where $C_{3/2}$ and $C_2$ stand for the appropriate Sobolev embedding constants (see inequality \eqref{cont_embedding} when $q \in \{\frac{3}{2}, 2\}$). Since $\varepsilon>0$ can be arbitrarily small, the previous inequalities immediately imply \eqref{vegtelenben}.

    On the account of \eqref{nullaban} and \eqref{vegtelenben}, we can conclude that 
    \begin{equation*}%\label{boundness}
	0 <\sup_{u \in W^{1,2}_{0,F}(\Omega)\setminus\{0\}}\frac{\mathcal{G}(u)}{\mathcal{N}(u)}<+\infty.
    \end{equation*}
    Therefore, one has
    $$0 < \lambda_2 \coloneqq
    \inf_{\substack{u\in W^{1,2}_{0,F}(\Omega) \\
	\mathcal G(u)> 0}}\frac{\mathcal N(u)}{2\mathcal G(u)} < +\infty.$$

    Now let us fix a parameter $\lambda > \lambda_2$. Then there exists a function $w_\lambda \in W^{1,2}_{0,F}(\Omega)$ such that $\mathcal{G}(w_\lambda)>0$ and
    $$\lambda > \frac{\mathcal{N}(w_\lambda)}{2\mathcal{G}(w_\lambda)} \geq \lambda_2.$$
    This yields that
    \begin{align*}
        C_\lambda^1 &\coloneqq \inf_{u \in W^{1,2}_{0,F}(\Omega)} \mathcal{E}_{\lambda}(u) \\
        &\leq  \mathcal{E}_{\lambda}(w_\lambda) = \frac{1}{2}\mathcal{N}(w_\lambda)-\lambda \mathcal{G}(w_\lambda)<0.
    \end{align*}
	
    Since the energy functional  $\mathcal{E}_{\lambda}$ is bounded from below and verifies the non-smooth Palais-Smale condition (see Sections \ref{Coercivity} \& \ref{PS}), it follows by Chang \cite[Theorem 3.5]{Chang} that $C_\lambda^1$ is a critical value of $\mathcal{E}_{\lambda}$, thus there exists a function $u_{\lambda}^1 \in W^{1,2}_{0,F}(\Omega)$ such that $C_\lambda^1 = \mathcal{E}_{\lambda}(u_{\lambda}^1)$ and $0 \in  \partial \mathcal{E}_{\lambda}(u_{\lambda}^1)$. In particular, we have $u_{\lambda}^1 \neq 0$, since $C_\lambda^1 = \mathcal{E}_{\lambda}(u_{\lambda}^1) < 0 = \mathcal{E}_{\lambda}(0)$, which also implies that $\lambda_2 > \lambda_1$.

    According to Section \ref{B}, $u_{\lambda}^1$ is a nontrivial nonnegative weak solution of $({\mathcal DI})_\lambda$ when $\lambda > \lambda_2 > \lambda_1$.

    \subsection{Second solution}

    This section provides another, minimax-type critical point for the energy functional $\mathcal{E}_{\lambda}$.
    
    Let $\lambda >\lambda_2$ and $\varepsilon>0$ be arbitrarily fixed. Estimate \eqref{G_estimate} yields that for every $p \in (2,2^*)$ and $u \in W^{1,2}_{0,F}(\Omega)$,
    \begin{align}\label{E_estimate}
        \mathcal{E}_\lambda(u) &= \frac{1}{2}\mathcal{N}(u) - \lambda \mathcal{G}(u) \nonumber\\
        &\geq  \left( \frac{1}{2} -\lambda\varepsilon C_2^2\right) \int_{\Omega} F^{*2}(x,Du(x))\mathrm{d}v_F  \\
        &-\lambda k_{\varepsilon} C_p^p  \left(\int_{\Omega} F^{*2}(x,Du(x))\mathrm{d}v_F\right)^\frac{p}{2}, \nonumber
    \end{align}
    where $C_2$ and $C_p$ denote the Sobolev embedding constants from \eqref{cont_embedding} for every  $p \in (2, 2^*)$.
    
    By choosing $\varepsilon$ small enough, namely $\varepsilon < \frac{1}{2\lambda C_2^2}$, we obtain that
    $$\rho_\lambda \coloneqq \left(\frac{\frac{1}{2}-\lambda\varepsilon C_2^2}{\lambda k_{\varepsilon}C_p^p}\right)^{\frac{1}{p-2}}>0.$$
	
    Since $\rho_{\lambda}>0$, relation \eqref{E_estimate} implies that
    $$\displaystyle\inf_{\substack{u\in W^{1,2}_{0,F}(\Omega) \\
    \mathcal{N}(u) = \rho_{\lambda}^2}}\mathcal{E}_{\lambda}(u) \geq \left(\frac{1}{2}-\lambda\varepsilon C_2^2\right)\rho_\lambda^2 > 0,$$		
    and we recall that $0 = \mathcal{E}_{\lambda}(0) >  \mathcal{E}_{\lambda}(w_\lambda)$.

    Accordingly, $\mathcal{E}_{\lambda}$ has the mountain pass geometry. Since $\mathcal{E}_{\lambda}$ also verifies the non-smooth Palais-Smale condition (see Section \ref{PS}), we can apply the non-smooth  mountain pass theorem (see e.g., Kristály, Motreanu and Varga \cite[Theorem 2]{KMV}), which implies that there exists a critical point $u_{\lambda}^2 \in W^{1,2}_{0,F}(\Omega)$, such that
    $0 \in \partial \mathcal E_{\lambda}(u_{\lambda}^2)$ and
    $$ C^2_\lambda \coloneqq  \mathcal E_{\lambda}(u_{\lambda}^2) = \inf_{\gamma \in \Gamma} \max_{t\in [0,1]}\mathcal{E}_{\lambda}(\gamma(t)),$$
    where 
    $$\Gamma = \bigg\{\gamma \in C\left([0,1];W^{1,2}_{0,F}(\Omega)\right): \gamma(0) = 0, \gamma(1) = w_\lambda \bigg\}.$$
	
    Considering the fact that 
    $$C^2_\lambda = \mathcal{E}_{\lambda}(u_{\lambda}^2)>0 = \mathcal{E}_{\lambda}(0) >  \mathcal{E}_{\lambda}(u_{\lambda}^1),$$
    we clearly have
    $u_{\lambda}^1 \neq u_{\lambda}^2 $ and $u_{\lambda}^2 \neq 0$. Taking into account Section \ref{B}, it follows that $u_{\lambda}^2$ is the second nontrivial nonnegative weak solution of the inclusion problem $({\mathcal DI})_\lambda$, which concludes the proof of $(\textbf{R}_2)$.

\section{Acknowledgment}

Á. Mester was supported by the ÚNKP-22-4 New National Excellence Program of the Ministry for Culture and Innovation from the source of the National Research, Development and Innovation Fund.

\end{document}